\documentclass[11pt]{article}
\usepackage[margin=1in]{geometry}
\usepackage{amsmath,amssymb,amsthm}
\usepackage{tikz}
\usetikzlibrary{arrows.meta}
\tikzset{vertex/.style={circle,draw,fill=white,minimum size=0.6cm,inner sep=0pt},
  block/.style={rectangle,draw,minimum width=1.1cm,minimum height=0.7cm}}
\newcommand{\avec}{\operatorname{avec}}
\newcommand{\diam}{\operatorname{diam}}
\newcommand{\ecc}{\operatorname{ecc}}
\theoremstyle{plain}
\newtheorem{theorem}{Theorem}
\newtheorem{lemma}{Lemma}[section]
\newtheorem{proposition}{Proposition}[section]
\newtheorem{conjecture}{Conjecture}[section]
\theoremstyle{remark}
\newtheorem{remark}{Remark}[section]

\begin{document}

\title{On the average eccentricity of Eulerian and bipartite Eulerian Digraphs}
\author{F.J. Osaye\\ Department of Mathematics and Statistics\\ Troy University, AL}
\date{}
\maketitle

\begin{abstract}
Let $D$ be an Eulerian digraph of order $n$ and minimum degree $\delta$.
We prove that the average eccentricity of $D$ satisfies
$\avec(D) \leq \frac{3n}{\delta+1}$, and that
$\avec(D) \leq \frac{2n}{\delta}+2$ if $D$ is bipartite.
We show that both bounds are best possible apart from an additive
constant and the $\delta+1$ versus $\delta+2$ gap (resp.\ the
$\delta$ versus $\delta+1$ gap in the bipartite case), where the
average eccentricity of the extremal examples shows that
$\avec = \frac{3n}{\delta+2}-2$ in general and
$\avec = \frac{2n}{\delta+1}-\frac{5}{2}$ in the bipartite case
for $\delta \geq 2$. We further conjecture that these values are
the true maxima up to an additive constant (for $\delta \geq 3$ in
the bipartite case).
\end{abstract}



\section{Introduction}

A digraph $D$ is \emph{Eulerian} if it is strong and
$\deg^+(v) = \deg^-(v)$ for every vertex $v$.
Eulerian digraphs form a large class containing all connected graphs
(with each edge replaced by two opposite arcs), and several bounds on
distance parameters that hold for undirected graphs are known to extend
to Eulerian digraphs but not to strong digraphs in general.
For example, the diameter bound for graphs of order $n$ and minimum
degree $\delta$,
\[
\diam(G) \leq \frac{3n}{\delta+1}+O(1),
\]
was shown by Soares \cite{Soares} and independently by Knyazev
\cite{Knyazev} to hold also for Eulerian digraphs.
More recently, it was shown that bounds on the size of graphs of
given order and diameter can also be extended to all Eulerian
digraphs \cite{DankelmannSize}.

However, the bound on the average eccentricity of an undirected
graph $G$ of order $n$ and minimum degree $\delta$ proven by
Dankelmann, Goddard and Swart \cite{DankelmannGoddardSwart} to be
$\avec(G) \leq \frac{9n}{4(\delta+1)}+O(1)$,
does not extend to Eulerian digraphs for $\delta \geq 3$.

Dankelmann \cite{Dankelmann} proved the corresponding bound on the
average distance $\mu(D)$ of an Eulerian digraph of order $n$ and
minimum degree $\delta$,
\[
\mu(D) \leq \frac{3n}{2(\delta+1)} + O(1),
\]
together with the improved bound $\mu(D) \leq \frac{n}{\delta}+O(1)$
for bipartite Eulerian digraphs.
His proof groups the distance layers from a vertex into blocks of three
(resp.\ four) and uses the Eulerian condition to show that every block
contains many vertices.

In this paper, we adapt the same method used by Dankelmann in
\cite{Dankelmann} to the average eccentricity.
Recall that the eccentricity of a vertex $v$ is
$\ecc(v) = \max_{x} d(v,x)$, where $d$ denotes directed distance,
and the average eccentricity is
$\avec(D) = \frac{1}{n}\sum_{v} \ecc(v)$.
While the average distance, diameter and the Wiener index of directed
graphs have been studied extensively in the literature, that of the
average eccentricity have received less attention. The first systematic
study of bounds on the average eccentricity of digraphs was carried out by
Osaye \cite{Osaye} who showed that for a $k$-connected tournament $T$,
\[
2 \leq \avec(T) < \frac{n}{2k} + \frac{11}{2} + \frac{k}{8n},
\]
and \[ 2\leq \avec(T) \leq \frac{n+11}{2} + \frac{1}{8n}\] if $T$ is a strong tournament. A known earlier work by Gimbert et. al \cite{GimbertLopezMillerRyan} on eccentricity in digraphs concerned the eccentric digraph ED(G)  in which  $u\longrightarrow v$ if and only if $v$ is at maximum distance from $u$ in $G$.
In particular, no bounds on the average eccentricity of digraphs were known before \cite{Osaye}. The purpose of this paper is to bridge this gap by proving further
upper bounds on the average eccentricity of Eulerian and bipartite
Eulerian digraphs. Our main results are the following.

\begin{theorem}
Let $D$ be an Eulerian digraph of order $n \geq 2$ and minimum degree
$\delta$. Then
\[
\avec(D) \leq \frac{3n}{\delta+1}.
\]
\end{theorem}

\begin{theorem}
Let $D$ be a bipartite Eulerian digraph of order $n \geq 2$ and
minimum degree $\delta$. Then
\[
\avec(D) \leq \frac{2n}{\delta}+2.
\]
\end{theorem}

Since $\avec(D) \leq \diam(D)$, Theorem 1 can also be read as a proof
of the diameter bound $\diam(D) < \frac{3n}{\delta+1}$ for Eulerian
digraphs by the layer method.
The point of the present note is that the same method, applied to
every vertex rather than to a diametral pair only, bounds the average
eccentricity by the same quantity.
In Section 4 we show that Theorem 1 is best possible apart from the
$\delta+1$ versus $\delta+2$ gap, and in Section 5 we do the same for
Theorem 2.

\section{Preliminaries}

All digraphs are finite and have no loops or multiple arcs.
For a digraph $D$, $V(D)$ denotes the vertex set and $n = |V(D)|$ the
order. For $v \in V(D)$, $\deg^+(v)$ and $\deg^-(v)$ denote the
out-degree and in-degree of $v$, and the minimum degree is
$\delta(D) = \min_{v} \min(\deg^+(v),\deg^-(v))$.
Since $D$ is Eulerian, $\deg^+(v) = \deg^-(v) \geq \delta$ for every $v$,
and $D$ strong implies $\delta \geq 1$ when $n \geq 2$.

For $v \in V(D)$ and $i \geq 0$, the $i$-th distance layer from $v$ is
\[
N_i(v) = \{x \in V(D) : d(v,x) = i\}, \qquad n_i(v) = |N_i(v)|.
\]
If $D$ is strong, then $N_i(v) \neq \emptyset$ if and only if
$0 \leq i \leq \ecc(v)$, and the layers $N_0(v),\dots,N_{\ecc(v)}(v)$
partition $V(D)$.
We write $N_{\leq i} = \bigcup_{j \leq i} N_j$ and
$N_{\geq i} = \bigcup_{j \geq i} N_j$, dropping the argument $v$ when
it is clear from context.
For $A,B \subseteq V(D)$, $m(A,B)$ denotes the number of arcs with tail
in $A$ and head in $B$.

We use the following consequence of the Eulerian condition throughout.

\begin{lemma}
Let $D$ be Eulerian and let $(A,B)$ be a partition of $V(D)$.
Then $m(A,B) = m(B,A)$.
\end{lemma}

\begin{proof}
Summing $\deg^+(v) - \deg^-(v) = 0$ over all $v \in A$ counts every arc
with both ends in $A$ once positively and once negatively, leaving
$m(A,B) - m(B,A) = 0$.
\end{proof}

\section{The general bound}

Fix a vertex $v$ and write $N_i = N_i(v)$, $n_i = n_i(v)$,
$k = \ecc(v)$. The following is the key estimate; it is Claim 1 in
the proof of Theorem 4 of \cite{Dankelmann}, proved here in full.

\begin{lemma}
\label{lem:triple}
For every $i \in \{0,1,\dots,k\}$,
\[
n_{i-1} + n_i + n_{i+1} \geq \delta+1,
\]
where $n_{-1} = n_{k+1} = 0$.
\end{lemma}

\begin{proof}
For $i = 0$, we have $N_0 = \{v\}$ and $N_1$ is exactly the set of
out-neighbours of $v$, so $n_1 = \deg^+(v) \geq \delta$ and
$n_{-1}+n_0+n_1 = 1+n_1 \geq \delta+1$.

Let $1 \leq i \leq k-1$.
Every out-neighbour $w$ of a vertex $u \in N_i$ satisfies
$d(v,w) \leq d(v,u)+1 = i+1$, so $w \in N_{\leq i+1}$.
Hence, summing out-degrees over $N_i$,
\[
\sum_{u \in N_i} \deg^+(u)
= m(N_i, N_{\leq i-1}) + m(N_i, N_i \cup N_{i+1}).
\]
For the second term, each $u \in N_i$ has at most
$n_i + n_{i+1} - 1$ out-neighbours in $N_i \cup N_{i+1}$ (all vertices
there except $u$ itself), so
$m(N_i, N_i \cup N_{i+1}) \leq n_i(n_i+n_{i+1}-1)$.
For the first term, $N_i \subseteq N_{\geq i}$, so
$m(N_i, N_{\leq i-1}) \leq m(N_{\geq i}, N_{\leq i-1})$.
By Lemma 2.1 applied to the partition
$(N_{\leq i-1}, N_{\geq i})$,
$m(N_{\geq i}, N_{\leq i-1}) = m(N_{\leq i-1}, N_{\geq i})$.
Moreover, every arc from $N_{\leq i-1}$ to $N_{\geq i}$ must go from
$N_{i-1}$ to $N_i$: if $x \in N_j$ with $j \leq i-2$ and $x \to y$ with
$y \in N_{\geq i}$, then $d(v,y) \leq j+1 \leq i-1$, a contradiction.
Thus $m(N_{\leq i-1}, N_{\geq i}) \leq n_{i-1} n_i$, and so
$m(N_i, N_{\leq i-1}) \leq n_{i-1} n_i$.
Combining,
\[
\delta n_i \leq \sum_{u \in N_i} \deg^+(u)
\leq n_{i-1} n_i + n_i(n_i+n_{i+1}-1).
\]
Since $n_i \geq 1$ (as $i \leq k$), dividing by $n_i$ gives
$\delta \leq n_{i-1}+n_i+n_{i+1}-1$, as claimed.

For $i = k$, we have $N_{\geq k} = N_k$ and
$\sum_{u \in N_k} \deg^+(u) = m(N_k, N_{\leq k-1}) + m(N_k, N_k)$.
As above, $m(N_k, N_{\leq k-1}) \leq m(N_{\leq k-1}, N_k) \leq
n_{k-1} n_k$ (arcs from $N_{\leq k-1}$ to $N_k$ must come from
$N_{k-1}$), and $m(N_k, N_k) \leq n_k(n_k-1)$ since there are no loops.
Thus $\delta n_k \leq n_{k-1} n_k + n_k(n_k-1)$, giving
$n_{k-1}+n_k \geq \delta+1$, which is the claim with $n_{k+1} = 0$.
\end{proof}

\begin{lemma}
\label{lem:ecc}
For every vertex $v$ of $D$,
\[
\ecc(v) < \frac{3n}{\delta+1}.
\]
\end{lemma}

\begin{proof}
Let $k = \ecc(v)$ and write $k+1 = 3q+r$ with $q \geq 0$ and
$r \in \{0,1,2\}$, so there are $q$ full triples of consecutive layers
and a remainder of $r$ layers.
For each $j$ with $3j+2 \leq k$, Lemma~\ref{lem:triple} applied at
$i = 3j+1$ gives
\[
n_{3j} + n_{3j+1} + n_{3j+2} \geq \delta+1.
\]
Hence the $q$ full triples contribute at least $q(\delta+1)$ vertices.
For the remainder: if $r = 0$ there is none; if $r = 1$ the remaining
layer is $N_k$ with $n_k \geq 1$; if $r = 2$ the remaining layers are
$N_{k-1}, N_k$ with $n_{k-1}+n_k \geq \delta+1$ by
Lemma~\ref{lem:triple} at $i = k$.
Therefore:
\begin{itemize}
\item If $r = 0$: $n \geq q(\delta+1)$, so
$k = 3q-1 \leq \frac{3n}{\delta+1}-1 < \frac{3n}{\delta+1}$.
\item If $r = 1$: $n \geq q(\delta+1)+1$, so
$k = 3q \leq \frac{3(n-1)}{\delta+1} < \frac{3n}{\delta+1}$.
\item If $r = 2$: $n \geq (q+1)(\delta+1)$, so
$k = 3q+1 = 3(q+1)-2 \leq \frac{3n}{\delta+1}-2
< \frac{3n}{\delta+1}$.
\end{itemize}
In all cases, $\ecc(v) = k < \frac{3n}{\delta+1}$.
\end{proof}

\begin{proof}[Proof of Theorem 1]
By Lemma~\ref{lem:ecc}, $\ecc(v) < \frac{3n}{\delta+1}$ for every
$v \in V(D)$. Averaging over all vertices,
\[
\avec(D) = \frac{1}{n}\sum_{v \in V(D)} \ecc(v)
< \frac{3n}{\delta+1},
\]
in particular $\avec(D) \leq \frac{3n}{\delta+1}$.
\end{proof}

\begin{remark}
Lemma~\ref{lem:ecc} also gives $\diam(D) < \frac{3n}{\delta+1}$,
since the diameter is the maximum eccentricity. Thus the layer method
recovers the Soares--Knyazev diameter bound for Eulerian digraphs.
\end{remark}

\section{Sharpness of Theorem 1}

We now show that Theorem 1 is best possible apart from the
$\delta+1$ versus $\delta+2$ gap, by computing the average eccentricity
of the extremal example from \cite{Dankelmann} exactly.

\begin{proposition}[Construction]
Let $\delta \geq 1$ and $k \geq 2$.
For $i = 1,\dots,k$, let $V_i = \{a_i,b_i,c_i\} \cup W_i$ with
$|W_i| = \delta-1$ (so $W_i = \emptyset$ if $\delta = 1$), and let
$D_i$ be the digraph on $V_i$ with all arcs $u \to v$ ($u \neq v$)
except the four arcs $a_i \to b_i$, $b_i \to c_i$, $c_i \to a_i$,
$b_i \to a_i$.
Let $D_{\delta,k}$ be obtained from the disjoint union of the $D_i$
by adding the arcs $b_i \to a_{i+1}$ for $i = 1,\dots,k$, with indices
taken modulo $k$.
Then $D_{\delta,k}$ is Eulerian, strong, of order $n = k(\delta+2)$,
and of minimum degree $\delta$.
\end{proposition}

\begin{proof}
Clearly $n = k(\delta+2)$.
In the complete digraph on $V_i$, every vertex has in-degree and
out-degree $\delta+1$.
The four removed arcs change the degrees as follows:
$a_i$ loses one out-arc ($a_i \to b_i$) and two in-arcs
($c_i \to a_i$, $b_i \to a_i$), so it has out-degree $\delta$ and
in-degree $\delta-1$ within $D_i$;
$b_i$ loses two out-arcs ($b_i \to c_i$, $b_i \to a_i$) and one in-arc
($a_i \to b_i$), so it has out-degree $\delta-1$ and in-degree
$\delta$ within $D_i$;
$c_i$ loses one out-arc ($c_i \to a_i$) and one in-arc ($b_i \to c_i$),
so it has out-degree and in-degree $\delta$ within $D_i$;
vertices of $W_i$ lose nothing.
Adding $b_i \to a_{i+1}$ increases the out-degree of $b_i$ by one and
the in-degree of $a_{i+1}$ by one.
Hence in $D_{\delta,k}$ every $a_i$, $b_i$, $c_i$ has in-degree and
out-degree exactly $\delta$, and every vertex of $W_i$ has in-degree
and out-degree $\delta+1$.
Thus $D_{\delta,k}$ is Eulerian with minimum degree $\delta$.

To see that $D_{\delta,k}$ is strong, note first that each $D_i$ is
strong when $\delta \geq 2$ (every missing arc is replaced by a path
of length $2$ through $W_i$, e.g.\ $a_i \to w \to b_i$), and that the
linking arcs $b_i \to a_{i+1}$ join the blocks in a directed cycle,
so the whole digraph is strong.
For $\delta = 1$, $D_i$ alone is not strong (for instance $b_i$ has
no out-neighbour in $V_i$), but the whole digraph still is: from
$a_i$ one reaches $c_i$ and $b_i$ within $V_i$ and then, via
$b_i \to a_{i+1}$, every other block; from $c_i$ one reaches $b_i$
and continues similarly; and from $b_i$ one goes to $a_{i+1}$ and
around the cycle back to $V_i$.
Hence $D_{\delta,k}$ is strong in all cases.
\end{proof}

Assume henceforth that $\delta \geq 2$, so that $W_i \neq \emptyset$;
the case $\delta = 1$ is discussed in Remark~\ref{rem:delta1}.
The key point is that the only arcs between distinct blocks are the
linking arcs $b_i \to a_{i+1}$.
Consequently, a shortest path never visits a block, leaves it, and
returns to it: leaving a block is possible only through its $b$-vertex,
and re-entering it is possible only through its $a$-vertex, so a path
that left and returned would contain a directed cycle and could be
shortened. We use this repeatedly below.

For $x \in V_i$, write $\alpha(x)$ for the distance from $x$ to $b_i$
within $D_i$, and for $y \in V_j$, write $\beta(y)$ for the distance
from $a_j$ to $y$ within $D_j$. Direct inspection of $D_i$ gives
\[
\begin{array}{c|cccc}
x & a_i & b_i & c_i & w \in W_i \\ \hline
\alpha(x) & 2 & 0 & 1 & 1
\end{array}
\qquad
\begin{array}{c|cccc}
y & a_j & b_j & c_j & w \in W_j \\ \hline
\beta(y) & 0 & 2 & 1 & 1
\end{array}
\]
Indeed, $a_i \to b_i$ is missing but $a_i \to c_i \to b_i$ has length
$2$; $b_i \to a_i$ and $b_i \to c_i$ are missing but
$b_i \to w \to a_i$ and $b_i \to w \to c_i$ have length $2$;
all other pairs are joined by an arc.
Moreover, $d(a_i,b_i) = 2$ within $D_i$, and the diameter of $D_i$
is $2$.

\begin{lemma}
\label{lem:dist2}
Let $x \in V_i$ and $y \in V_j$, and let $t \in \{0,\dots,k-1\}$ be
such that $j \equiv i+t \pmod{k}$.
If $t \geq 1$, then
\[
d(x,y) = \alpha(x) + \beta(y) + 3t - 2.
\]
If $t = 0$ (so $x,y$ lie in the same block), then $d(x,y) \leq 2$.
\end{lemma}

\begin{proof}
If $t = 0$, both vertices lie in $V_i = V_j$, and the diameter of
$D_i$ is $2$, as noted above.

Let $t \geq 1$.
We first exhibit a path of the claimed length:
$x \to b_i$ in $\alpha(x)$ steps within $D_i$,
then $b_i \to a_{i+1}$ (one step),
then $a_{i+1} \to b_{i+1}$ in $2$ steps within $D_{i+1}$,
and so on through the blocks, ending with
$b_{j-1} \to a_j$ (one step) and $a_j \to y$ in $\beta(y)$ steps
within $D_j$.
There are $t$ linking arcs and $t-1$ intermediate transits of length
$2$, for a total of $\alpha(x) + t + 2(t-1) + \beta(y)
= \alpha(x) + \beta(y) + 3t - 2$.

Conversely, let $P$ be a shortest $x \to y$ path.
Since the only arcs between distinct blocks are the linking arcs,
$P$ must use at least $t$ of them, in the order
$b_i \to a_{i+1}, \dots$ (possibly with extra wraps around the cycle);
write $L \geq t$ for the number used.
The initial segment of $P$, from $x$ to $b_i$, stays in $V_i$ by the
remark above (the only exit from $V_i$ is $b_i$ itself), so it has
length at least $\alpha(x)$.
Similarly, the final segment, from $a_j$ to $y$, stays in $V_j$ and
has length at least $\beta(y)$.
Between consecutive linking arcs, $P$ goes from $a_\ell$ to $b_\ell$
for some $\ell$; this segment stays in $V_\ell$ (leaving would force
a return through $a_\ell$) and has length at least
$d(a_\ell,b_\ell) = 2$.
There are $L-1$ such segments.
Hence $|P| \geq \alpha(x) + L + 2(L-1) + \beta(y)
= \alpha(x) + \beta(y) + 3L - 2
\geq \alpha(x) + \beta(y) + 3t - 2$.
Combined with the construction, this proves the formula.
\end{proof}

\begin{proposition}
\label{prop:avec2}
For $\delta \geq 2$ and $k \geq 2$,
\[
\diam(D_{\delta,k}) = 3k-1, \qquad
\avec(D_{\delta,k}) = 3k-2 = \frac{3n}{\delta+2}-2,
\]
where $n = k(\delta+2)$.
\end{proposition}

\begin{proof}
Let $x \in V_i$.
By Lemma~\ref{lem:dist2}, distances to vertices in the same block are
at most $2$, while distances to block $j$ with forward distance $t$
are $\alpha(x)+\beta(y)+3t-2$, which is maximized for $t = k-1$ and
$y = b_{i-1}$ (where $\beta(y) = 2$ is largest).
Since $3k-3 \geq 3 > 2$ for $k \geq 2$, the maximum is attained at
$t = k-1$, giving
\[
\ecc(x) = \alpha(x) + 2 + 3(k-1) - 2 = \alpha(x) + 3k - 3.
\]
Thus $\ecc(a_i) = 3k-1$, $\ecc(b_i) = 3k-3$,
$\ecc(c_i) = 3k-2$, and $\ecc(w) = 3k-2$ for $w \in W_i$.
In particular $\diam(D_{\delta,k}) = 3k-1$.
Averaging over one block,
\[
\frac{(3k-1)+(3k-3)+(3k-2)+(\delta-1)(3k-2)}{\delta+2}
= \frac{2(3k-2)+\delta(3k-2)}{\delta+2} = 3k-2,
\]
since $(3k-1)+(3k-3) = 2(3k-2)$.
All blocks are symmetric, so $\avec(D_{\delta,k}) = 3k-2$.
Finally, $n = k(\delta+2)$ gives $3k-2 = \frac{3n}{\delta+2}-2$.
\end{proof}

\begin{remark}
\label{rem:delta1}
For $\delta = 1$, the sets $W_i$ are empty and the within-block
shortcuts through $W_i$ disappear; a direct check using
Lemma~\ref{lem:dist2} (with same-block targets reached by going
around the cycle) shows that every vertex of $D_{1,k}$ has
eccentricity $3k-1$, so $\avec(D_{1,k}) = 3k-1
= \frac{3n}{\delta+2}-1$.
In particular, $\avec(D_{\delta,k}) = \frac{3n}{\delta+2}+O(1)$ for
all $\delta \geq 1$.
\end{remark}

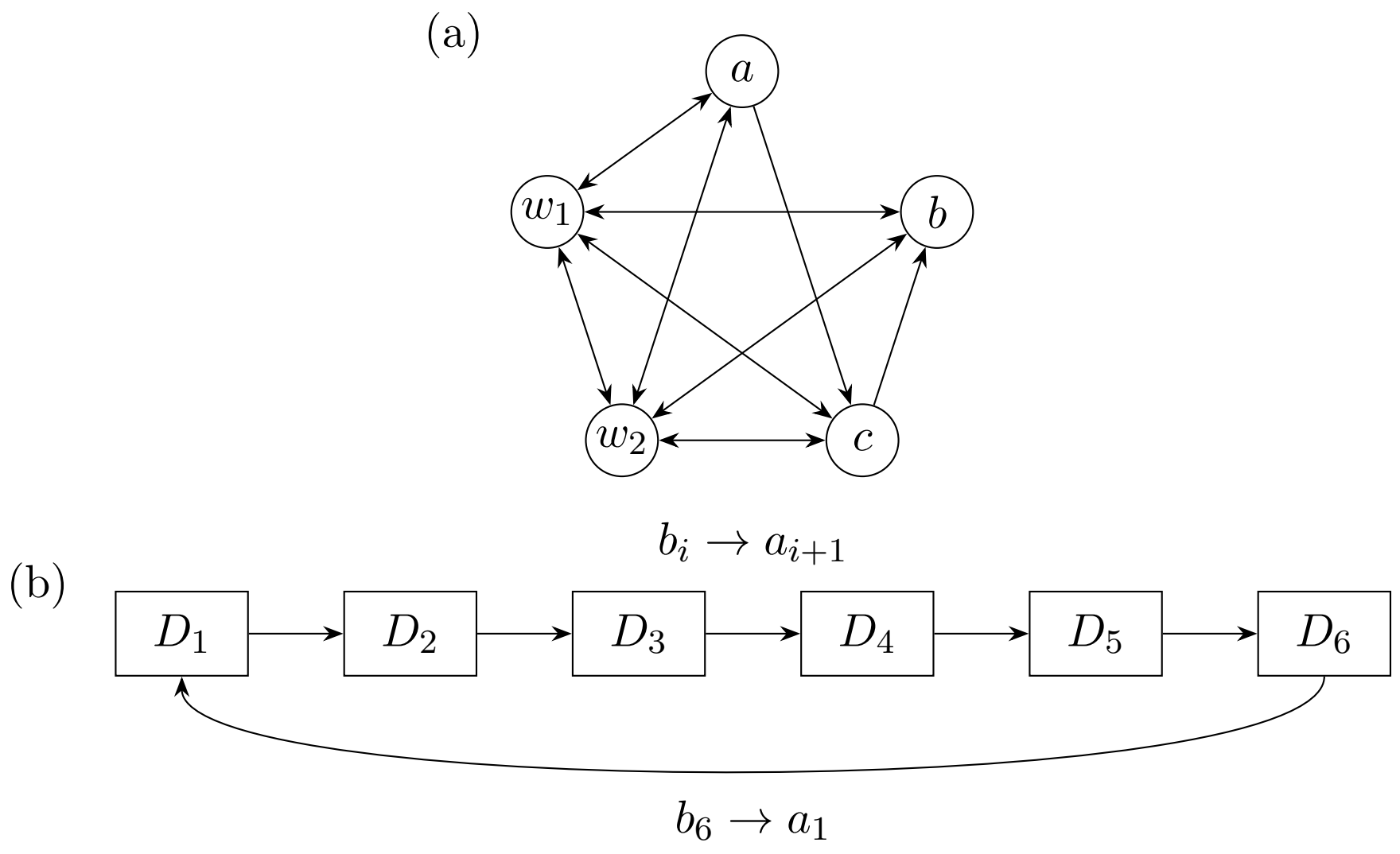
\begin{figure}[ht]
\centering
\begin{tikzpicture}[>=Stealth]
  \node at (-2.4,2.0) {(a)};
  \node[vertex] (a) at (0,1.7) {$a$};
  \node[vertex] (b) at (1.62,0.53) {$b$};
  \node[vertex] (c) at (1.0,-1.37) {$c$};
  \node[vertex] (w2) at (-1.0,-1.37) {$w_2$};
  \node[vertex] (w1) at (-1.62,0.53) {$w_1$};
  \draw[<->] (a) -- (w1);
  \draw[<->] (a) -- (w2);
  \draw[<->] (b) -- (w1);
  \draw[<->] (b) -- (w2);
  \draw[<->] (c) -- (w1);
  \draw[<->] (c) -- (w2);
  \draw[<->] (w1) -- (w2);
  \draw[->] (a) -- (c);
  \draw[->] (c) -- (b);
\end{tikzpicture}
\\[6pt]
\begin{tikzpicture}[>=Stealth]
  \node at (-1.2,0.4) {(b)};
  \node[block] (D1) at (0,0) {$D_1$};
  \node[block] (D2) at (1.9,0) {$D_2$};
  \node[block] (D3) at (3.8,0) {$D_3$};
  \node[block] (D4) at (5.7,0) {$D_4$};
  \node[block] (D5) at (7.6,0) {$D_5$};
  \node[block] (D6) at (9.5,0) {$D_6$};
  \draw[->] (D1) -- (D2);
  \draw[->] (D2) -- (D3);
  \draw[->] (D3) -- (D4);
  \draw[->] (D4) -- (D5);
  \draw[->] (D5) -- (D6);
  \draw[->] (D6.south) .. controls (9.5,-1.4) and (0,-1.4) .. (D1.south);
  \node at (4.75,0.75) {$b_i \to a_{i+1}$};
  \node at (4.75,-1.55) {$b_6 \to a_1$};
\end{tikzpicture}
\caption{The extremal digraph $D_{3,6}$: (a) one block $D_i$ for
$\delta = 3$ ($W_i = \{w_1,w_2\}$; double-headed arrows denote arcs
in both directions; the four missing arcs are $a \to b$, $b \to a$,
$b \to c$, $c \to a$); (b) the six blocks joined by the linking
arcs $b_i \to a_{i+1}$.}
\label{fig:D36}
\end{figure}

Figure~\ref{fig:D36} shows the extremal digraph $D_{3,6}$:
one block $D_i$ in full, together with the six-block linking cycle.
For $D_{3,6}$ we have $n = 30$, $\delta = 3$, $\diam = 17$ and
$\avec = 16 = \frac{3n}{\delta+2}-2$.

The bound of Theorem 1 is therefore best possible apart from the
$\delta+1$ versus $\delta+2$ gap and an additive constant.
This suggests the following analogue of Dankelmann's Conjecture 1
for the average distance.

\begin{conjecture}
There exists a constant $c$ such that every Eulerian digraph $D$ of
order $n$ and minimum degree $\delta$ satisfies
\[
\avec(D) \leq \frac{3n}{\delta+2} + c.
\]
\end{conjecture}

\section{The bipartite case}

Let now $D$ be bipartite and Eulerian: the underlying undirected graph
is bipartite with parts $X$ and $Y$, so every arc joins $X$ to $Y$ or
$Y$ to $X$.
For $v \in V(D)$, the distance layers alternate between the parts:
$N_i(v) \subseteq X$ for $i$ even and $N_i(v) \subseteq Y$ for $i$ odd
(or vice versa, depending on the part containing $v$).
In particular, there are no arcs within a single layer $N_i$.

\begin{lemma}
\label{lem:bip}
Let $v \in V(D)$, $k = \ecc(v)$, $n_i = |N_i(v)|$.
Then for every $i \in \{0,1,\dots,k\}$,
\[
n_{i-1} + n_{i+1} \geq \delta,
\]
where $n_{-1} = n_{k+1} = 0$.
\end{lemma}

\begin{proof}
For $i = 0$, $n_1 = \deg^+(v) \geq \delta$, so
$n_{-1}+n_1 \geq \delta$.

Let $1 \leq i \leq k-1$.
Every out-neighbour of $u \in N_i$ lies in $N_{i+1} \cup N_{\leq i-1}$,
and, since there are no arcs within $N_i$,
\[
\sum_{u \in N_i} \deg^+(u)
= m(N_i, N_{\leq i-1}) + m(N_i, N_{i+1}).
\]
As in Lemma~\ref{lem:triple},
$m(N_i, N_{\leq i-1}) \leq m(N_{\geq i}, N_{\leq i-1})
= m(N_{\leq i-1}, N_{\geq i}) \leq n_{i-1} n_i$,
using Lemma 2.1 and the fact that arcs from $N_{\leq i-1}$ to
$N_{\geq i}$ must come from $N_{i-1}$.
Also $m(N_i, N_{i+1}) \leq n_i n_{i+1}$.
Thus $\delta n_i \leq n_{i-1} n_i + n_i n_{i+1}$, and dividing by
$n_i \geq 1$ gives $n_{i-1}+n_{i+1} \geq \delta$.

For $i = k$, we have
$\sum_{u \in N_k} \deg^+(u) = m(N_k, N_{\leq k-1}) + m(N_k, N_k)$,
but $m(N_k, N_k) = 0$ since there are no arcs within a layer.
As before, $m(N_k, N_{\leq k-1}) \leq n_{k-1} n_k$, so
$\delta n_k \leq n_{k-1} n_k$, giving $n_{k-1} \geq \delta$,
which is the claim with $n_{k+1} = 0$.
\end{proof}

\begin{lemma}
\label{lem:bipecc}
For every vertex $v$ of a bipartite Eulerian digraph $D$,
\[
\ecc(v) < \frac{2n}{\delta}+2.
\]
\end{lemma}

\begin{proof}
Let $k = \ecc(v)$ and write $k+1 = 4q+r$ with $q \geq 0$ and
$r \in \{0,1,2,3\}$, so there are $q$ full blocks of four consecutive
layers and a remainder of $r$ layers.
For each full block $\{4j,4j+1,4j+2,4j+3\}$ (so $4j+3 \leq k$),
Lemma~\ref{lem:bip} at $i = 4j+1$ gives
$n_{4j}+n_{4j+2} \geq \delta$, and at $i = 4j+2$ gives
$n_{4j+1}+n_{4j+3} \geq \delta$; hence the block contains at least
$2\delta$ vertices.
Each of the $r$ remaining layers contains at least one vertex.
Thus $n \geq 2\delta q + r$, so
$4q \leq \frac{2(n-r)}{\delta}$ and
\[
k = 4q+r-1 \leq \frac{2n}{\delta} + r\left(1-\frac{2}{\delta}\right) - 1.
\]
The right-hand side is maximized over $r \in \{0,1,2,3\}$ by $r = 3$,
giving $k < \frac{2n}{\delta}+2$ (indeed, for $r = 3$ the bound is
$\frac{2n}{\delta}+2-\frac{6}{\delta} < \frac{2n}{\delta}+2$,
and the other cases give even smaller values).
Hence $\ecc(v) < \frac{2n}{\delta}+2$.
\end{proof}

\begin{proof}[Proof of Theorem 2]
By Lemma~\ref{lem:bipecc}, $\ecc(v) < \frac{2n}{\delta}+2$ for every
$v \in V(D)$. Averaging gives
$\avec(D) \leq \frac{2n}{\delta}+2$.
\end{proof}

We now show that Theorem 2 is best possible apart from the
$\delta$ versus $\delta+1$ gap, using the bipartite extremal example
from \cite{Dankelmann}.

\begin{proposition}[Construction]
Let $\delta \geq 1$ and $k \geq 2$.
For $i = 1,\dots,k$, let $X_i = \{a_i,c_i\} \cup X'_i$ and
$Y_i = \{b_i,d_i\} \cup Y'_i$ with $|X'_i| = |Y'_i| = \delta-1$,
and let $D'_i$ be the digraph on $X_i \cup Y_i$ with all arcs from
$X_i$ to $Y_i$ and from $Y_i$ to $X_i$ except the five arcs
$a_i \to b_i$, $b_i \to c_i$, $c_i \to d_i$, $d_i \to a_i$,
$b_i \to a_i$.
Let $D'_{\delta,k}$ be obtained from the disjoint union of the
$D'_i$ by adding the arcs $b_i \to a_{i+1}$ for $i = 1,\dots,k$,
with indices taken modulo $k$.
Then $D'_{\delta,k}$ is bipartite and Eulerian, strong, of order
$n = k(2\delta+2)$, and of minimum degree $\delta$.
\end{proposition}

\begin{proof}
The bipartition is $X = \bigcup_i X_i$, $Y = \bigcup_i Y_i$:
all arcs within blocks join $X_i$ to $Y_i$, and the linking arcs
$b_i \to a_{i+1}$ join $Y$ to $X$.
Clearly $n = k \cdot 2(\delta+1) = k(2\delta+2)$.
In the complete bipartite digraph on $(X_i,Y_i)$, every vertex has
in-degree and out-degree $\delta+1$.
The five removed arcs change the degrees as follows:
$a_i$ loses one out-arc ($a_i \to b_i$) and two in-arcs
($d_i \to a_i$, $b_i \to a_i$);
$b_i$ loses two out-arcs ($b_i \to c_i$, $b_i \to a_i$) and one in-arc
($a_i \to b_i$);
$c_i$ loses one out-arc ($c_i \to d_i$) and one in-arc ($b_i \to c_i$);
$d_i$ loses one out-arc ($d_i \to a_i$) and one in-arc ($c_i \to d_i$);
vertices of $X'_i \cup Y'_i$ lose nothing.
Adding $b_i \to a_{i+1}$ restores one out-degree to $b_i$ and one
in-degree to $a_{i+1}$.
Hence every $a_i$, $b_i$, $c_i$, $d_i$ has in-degree and out-degree
exactly $\delta$, and every vertex of $X'_i \cup Y'_i$ has in-degree
and out-degree $\delta+1$.
Thus $D'_{\delta,k}$ is Eulerian with minimum degree $\delta$.
Strong connectivity follows as in Proposition 4.1: when
$\delta \geq 2$ each block $D'_i$ has diameter $3$, hence is strong,
and the linking arcs join the blocks in a directed cycle;
for $\delta = 1$, every vertex reaches $b_i$ (or already is $b_i$)
and then continues around the linking cycle, returning to its block
via $a_i$, so the whole digraph is strong.
\end{proof}

Assume henceforth that $\delta \geq 2$.
As in Section 4, the only arcs between distinct blocks are the linking
arcs, so a shortest path never leaves a block and returns to it.
For $x \in X_i \cup Y_i$, write $\alpha(x)$ for the distance from $x$
to $b_i$ within $D'_i$, and for $y \in X_j \cup Y_j$, write
$\beta(y)$ for the distance from $a_j$ to $y$ within $D'_j$.
Direct inspection gives
\[
\begin{array}{c|cccccc}
x & a_i & b_i & c_i & d_i & x' \in X'_i & y' \in Y'_i \\ \hline
\alpha(x) & 3 & 0 & 1 & 2 & 1 & 2
\end{array}
\qquad
\begin{array}{c|cccccc}
y & a_j & b_j & c_j & d_j & x' \in X'_j & y' \in Y'_j \\ \hline
\beta(y) & 0 & 3 & 2 & 1 & 2 & 1
\end{array}
\]
For instance, $d(a_i,b_i) = 3$: the arc $a_i \to b_i$ is missing, a
path of length $2$ from $X$ to $Y$ is impossible in a bipartite
digraph (it would need a $Y \to Y$ or $X \to X$ step), and
$a_i \to d_i \to c_i \to b_i$ has length $3$.
Similarly $d(a_i,b_i) = 3$ is the transit cost between consecutive
linking arcs, and the diameter of $D'_i$ is $3$.

\begin{lemma}
\label{lem:dist3}
Let $x \in X_i \cup Y_i$ and $y \in X_j \cup Y_j$, and let
$t \in \{0,\dots,k-1\}$ be such that $j \equiv i+t \pmod{k}$.
If $t \geq 1$, then
\[
d(x,y) = \alpha(x) + \beta(y) + 4t - 3.
\]
If $t = 0$, then $d(x,y) \leq 3$.
\end{lemma}

\begin{proof}
If $t = 0$, both vertices lie in the same block, whose diameter is
$3$ as noted above.

Let $t \geq 1$.
The path $x \to b_i$ ($\alpha(x)$ steps), then $b_i \to a_{i+1}$,
then $a_{i+1} \to b_{i+1}$ ($3$ steps), and so on, ending with
$b_{j-1} \to a_j$ and $a_j \to y$ ($\beta(y)$ steps), has length
$\alpha(x) + t + 3(t-1) + \beta(y)
= \alpha(x) + \beta(y) + 4t - 3$.
Conversely, a shortest $x \to y$ path uses at least $t$ linking arcs
(extra wraps around the cycle would only lengthen it); write
$L \geq t$ for the number used.
Its initial segment stays in block $i$ (the only exit is $b_i$
itself) and has length at least $\alpha(x)$; its final segment stays
in block $j$ and has length at least $\beta(y)$; and between
consecutive linking arcs it goes from $a_\ell$ to $b_\ell$ within
block $\ell$, which costs at least $3$.
There are $L-1$ such segments, so its length is at least
$\alpha(x) + L + 3(L-1) + \beta(y)
= \alpha(x) + \beta(y) + 4L - 3
\geq \alpha(x) + \beta(y) + 4t - 3$,
proving the formula.
\end{proof}

\begin{proposition}
\label{prop:avec3}
For $\delta \geq 2$ and $k \geq 2$,
\[
\diam(D'_{\delta,k}) = 4k-1, \qquad
\avec(D'_{\delta,k}) = 4k-\frac{5}{2}
= \frac{2n}{\delta+1}-\frac{5}{2},
\]
where $n = k(2\delta+2)$.
\end{proposition}

\begin{proof}
Let $x$ lie in block $i$.
By Lemma~\ref{lem:dist3}, distances within the block are at most $3$,
while distances to block $j$ with forward distance $t$ are
$\alpha(x)+\beta(y)+4t-3$, maximized for $t = k-1$ and $y = b_{i-1}$
(where $\beta(y) = 3$ is largest).
Since $4k-4 \geq 4 > 3$ for $k \geq 2$, the maximum is attained at
$t = k-1$, giving
\[
\ecc(x) = \alpha(x) + 3 + 4(k-1) - 3 = \alpha(x) + 4k - 4.
\]
Thus $\ecc(a_i) = 4k-1$, $\ecc(b_i) = 4k-4$, $\ecc(c_i) = 4k-3$,
$\ecc(d_i) = 4k-2$, $\ecc(x') = 4k-3$ for $x' \in X'_i$, and
$\ecc(y') = 4k-2$ for $y' \in Y'_i$.
In particular $\diam(D'_{\delta,k}) = 4k-1$.
Averaging over one block,
\[
\begin{aligned}
&\frac{(4k-1)+(4k-4)+(4k-3)+(4k-2)
+(\delta-1)(4k-3)+(\delta-1)(4k-2)}{2(\delta+1)} \\
&\qquad = \frac{(\delta+1)(8k-5)}{2(\delta+1)} = 4k-\frac{5}{2},
\end{aligned}
\]
since $(4k-1)+(4k-4) = 8k-5$ and $(4k-3)+(4k-2) = 8k-5$.
All blocks are symmetric, so $\avec(D'_{\delta,k}) = 4k-\frac{5}{2}$.
Finally, $n = k(2\delta+2)$ gives
$4k-\frac{5}{2} = \frac{2n}{\delta+1}-\frac{5}{2}$.
\end{proof}

\begin{remark}
For $\delta = 1$, the sets $X'_i$ and $Y'_i$ are empty and a direct
check analogous to Remark~\ref{rem:delta1} shows that every vertex
of $D'_{1,k}$ has eccentricity $4k-1$, so
$\avec(D'_{1,k}) = 4k-1 = \frac{2n}{\delta+1}-1$.
In particular, $\avec(D'_{\delta,k}) = \frac{2n}{\delta+1}+O(1)$ for
all $\delta \geq 1$.
\end{remark}

\begin{figure}[ht]
\centering
\begin{tikzpicture}[>=Stealth]
  \node at (-1.0,3.4) {(a)};
  \node[vertex] (a) at (0,3) {$a$};
  \node[vertex] (c) at (0,2) {$c$};
  \node[vertex] (x1) at (0,1) {$x_1$};
  \node[vertex] (x2) at (0,0) {$x_2$};
  \node[vertex] (b) at (4.5,3) {$b$};
  \node[vertex] (d) at (4.5,2) {$d$};
  \node[vertex] (y1) at (4.5,1) {$y_1$};
  \node[vertex] (y2) at (4.5,0) {$y_2$};
  \draw[->] (a) -- (d);
  \draw[<->] (a) -- (y1);
  \draw[<->] (a) -- (y2);
  \draw[->] (c) -- (b);
  \draw[->] (d) -- (c);
  \draw[<->] (c) -- (y1);
  \draw[<->] (c) -- (y2);
  \draw[<->] (x1) -- (b);
  \draw[<->] (x1) -- (d);
  \draw[<->] (x1) -- (y1);
  \draw[<->] (x1) -- (y2);
  \draw[<->] (x2) -- (b);
  \draw[<->] (x2) -- (d);
  \draw[<->] (x2) -- (y1);
  \draw[<->] (x2) -- (y2);
\end{tikzpicture}
\\[7pt]
\begin{tikzpicture}[>=Stealth]
  \node at (-1.2,0.4) {(b)};
  \node[block] (D1) at (0,0) {$D'_1$};
  \node[block] (D2) at (1.9,0) {$D'_2$};
  \node[block] (D3) at (3.8,0) {$D'_3$};
  \node[block] (D4) at (5.7,0) {$D'_4$};
  \draw[->] (D1) -- (D2);
  \draw[->] (D2) -- (D3);
  \draw[->] (D3) -- (D4);
  \draw[->] (D4.south) .. controls (5.7,-1.4) and (0,-1.4) .. (D1.south);
  \node at (2.85,0.75) {$b_i \to a_{i+1}$};
  \node at (2.85,-1.55) {$b_4 \to a_1$};
\end{tikzpicture}
\caption{The bipartite extremal digraph $D'_{3,4}$: (a) one block
$D'_i$ for $\delta = 3$ ($X'_i = \{x_1,x_2\}$,
$Y'_i = \{y_1,y_2\}$; double-headed arrows denote arcs in both
directions; the five missing arcs are $a \to b$, $b \to a$,
$b \to c$, $c \to d$, $d \to a$); (b) the four blocks joined by the
linking arcs $b_i \to a_{i+1}$.}
\label{fig:Dp34}
\end{figure}
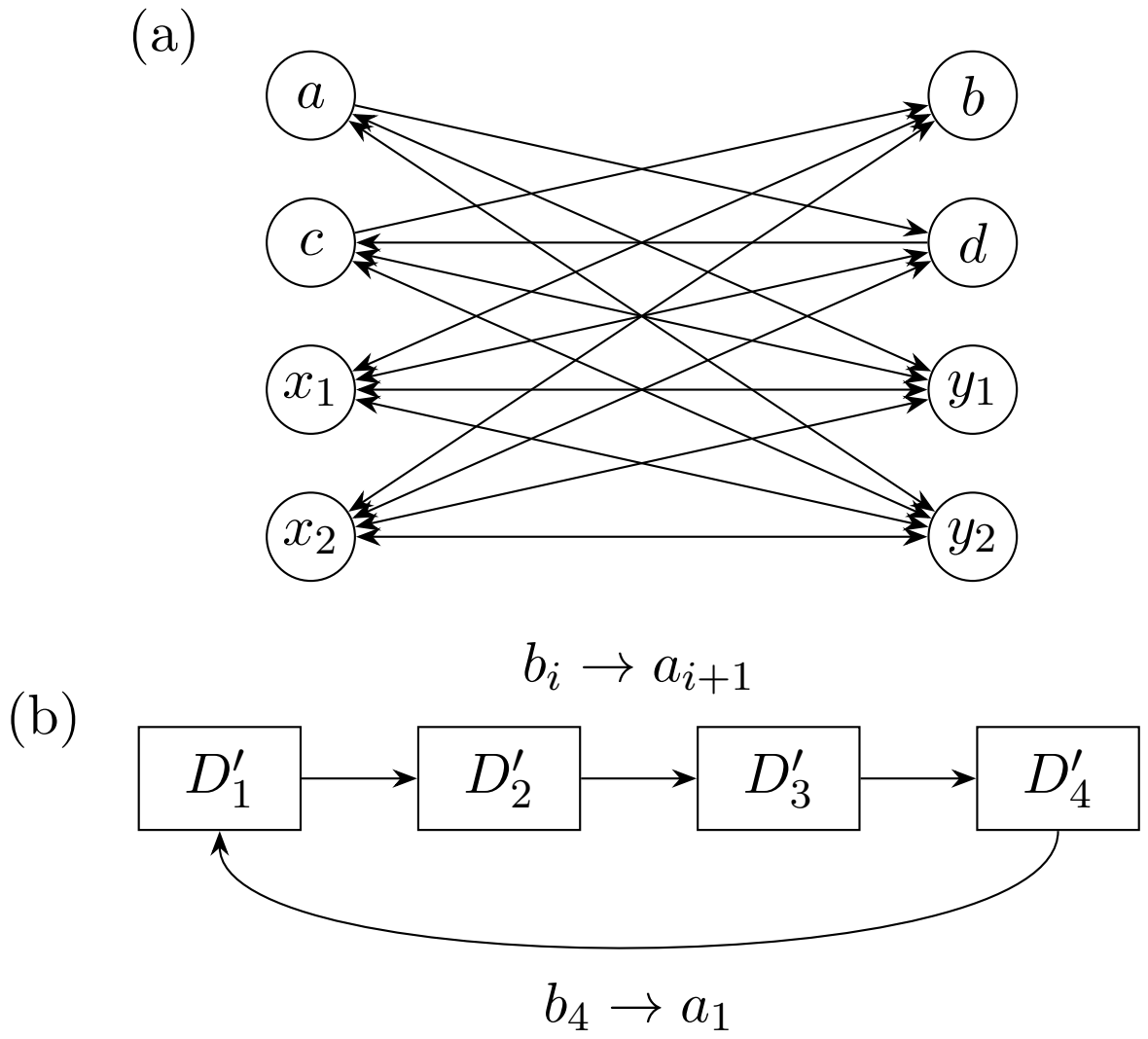

Figure~\ref{fig:Dp34} shows the bipartite extremal digraph
$D'_{3,4}$: one block $D'_i$ in full, together with the four-block
linking cycle. For $D'_{3,4}$ we have $n = 32$, $\delta = 3$,
$\diam = 15$ and $\avec = \frac{27}{2}
= \frac{2n}{\delta+1}-\frac{5}{2}$.

\begin{conjecture}
There exists a constant $c$ such that every bipartite Eulerian
digraph $D$ of order $n$ and minimum degree $\delta \geq 3$ satisfies
\[
\avec(D) \leq \frac{2n}{\delta+1} + c.
\]
\end{conjecture}

\begin{remark}
The restriction $\delta \geq 3$ is necessary: for $\delta = 2$ the
bound $\frac{2n}{3} + O(1)$ is false. Indeed, let blocks
$B_0, B_1, \dots, B_k$ have sizes $|B_0| = 1$, $|B_1| = 2$,
$|B_j| = 1$ for $2 \leq j \leq k-2$, and $|B_{k-1}| = |B_k| = 2$,
and join consecutive blocks by all arcs in both directions, with no
arcs within a block. The resulting digraph is bipartite, Eulerian and
has minimum degree $2$. The distance between vertices in blocks $B_j$
and $B_{j'}$ is $|j-j'|$ (and $2$ within a block of size $2$), so a
vertex in $B_j$ has eccentricity $\max(j, k-j)$, and averaging over
the blocks gives $\avec = \frac{3k}{4} + O(1) = \frac{3n}{4} + O(1)$,
which exceeds $\frac{2n}{3} + c$ for large $n$. Hence for
$\delta = 2$ the sharp bipartite bound is $\frac{3n}{4} + O(1)$,
the same as in the general case. (For $\delta = 1$ the bound
$\frac{2n}{\delta+1} + c$ holds trivially, since
$\avec \leq n-1$.)
\end{remark}

\section{Concluding remarks}

The bounds of Theorems 1 and 2 are obtained from the distance-layer
method only; no additional hypothesis on $D$ is used.
The remaining open problem is to close the $\delta+1$ versus
$\delta+2$ gap in Theorem 1 (resp.\ the $\delta$ versus $\delta+1$
gap in Theorem 2), as conjectured above. Note that, unlike for
tournaments, the average eccentricity of the extremal Eulerian
digraphs is essentially equal to the diameter bound, not half of it:
the extremal digraphs are cycle-like, so almost every vertex is
nearly diametral.

\section*{Acknowledgement}
The author would like to sincerely thank Professor Peter Dankelmann for suggesting research on the average eccentricity of digraphs.


\begin{thebibliography}{99}

\bibitem{Dankelmann} P. Dankelmann, On average distance in tournaments
and Eulerian digraphs, \emph{Discrete Appl.\ Math.} 266 (2019), 38--47.

\bibitem{DankelmannSize} P. Dankelmann, Distance and size in digraphs,
\emph{Discrete Math.} 338 (2015), 144--148.

\bibitem{DankelmannGoddardSwart} P. Dankelmann, W. Goddard, C.S. Swart,
The average eccentricity of a graph and its subgraphs.
\emph{Util.\ Math.} 41 (2004), 41--51.

\bibitem{GimbertLopezMillerRyan} J. Gimbert, N. L\'opez, M. Miller,
J. Ryan, Characterization of eccentric digraphs,
\emph{Discrete Math.} 306 (2006), 210--219.

\bibitem{Knyazev} A. V. Knyazev, Diameters of pseudosymmetric graphs,
\emph{Math.\ Notes} 41 (1987), 473--482.

\bibitem{Osaye} F. J. Osaye, Average eccentricity in $k$-connected tournaments (2026).
https://arxiv.org/pdf/2610.35970v1.

\bibitem{Soares} J. Soares, Maximum diameter of regular digraphs,
\emph{J.\ Graph Theory} 16 (1992), 437--450.






\end{thebibliography}
\end{document}